\documentclass[11pt]{scrartcl}

\usepackage{algorithm}
\usepackage{algorithmic}
\usepackage{amssymb}
\usepackage{amsmath}
\usepackage{amsthm}

\newcommand{\e}{\epsilon}
\newcommand{\be}{\unrhd_b}
\newcommand{\se}{\succ_\e}
\newcommand{\up}{\mathsf{update}}
\newcommand{\A}{\mathcal{A}}
\newcommand{\B}{\mathcal{B}}
\newcommand{\E}{\mathcal{E}}
\newcommand{\F}{\mathcal{F}}
\newcommand{\X}{\mathcal{X}}
\newcommand{\Xh}{\hat{\mathcal{X}}}
\newcommand{\Y}{\mathcal{Y}}
\newcommand{\Prob}{\mathbb{P}}
\newcommand{\T}{\mathsf{T}}
\newcommand{\Nat}{\mathbb{N}}
\newcommand{\Real}{\mathbb{R}}
\newcommand{\Zed}{\mathbb{Z}}
\newcommand{\lra}{\longrightarrow}

\newtheorem{definition}{Definition}
\newtheorem{lemma}{Lemma}
\newtheorem{theorem}{Theorem}

\begin{document}

\title{Stochastic convergence of random search to fixed size Pareto
  set approximations}

\author{Marco Laumanns\footnote{Institute for Operations Research,
    ETH Zurich, 8092 Zurich, Switzerland,
    \texttt{laumanns@ifor.math.ethz.ch}}} 

\maketitle

\begin{abstract}
This paper presents the first convergence result for random search
algorithms to a subset of the Pareto set of given maximum size $k$ with
bounds on the approximation quality $\e$.
The core of the algorithm is a new selection criterion based on a
hypothetical multilevel grid on the objective space.
It is shown that, when using this criterion for accepting new search
points, the sequence of solution archives converges with probability
one to a subset of the Pareto set that $\e$-dominates the entire
Pareto set. The obtained approximation quality $\e$ is equal to the
size of the grid cells on the finest level of resolution that allows an
approximation with at most $k$ points in the family of grids
considered. 
While the convergence result is of general theoretical interest, the
archiving algorithm might be of high practical value for any type
iterative multiobjective optimization method, such as evolutionary
algorithms or other metaheuristics, which all rely on the usage of a 
finite on-line memory to store the best solutions found so far as the
current approximation of the Pareto set.
\end{abstract}

\section{Introduction}

In multiobjective optimization, we are given $m \ge 2$ objective
functions $f_i : \X \lra \Real$, $i \in \{1,\ldots,m\}$, defined over
some search space $\X$, which might be implicitly defined by 
constraints. We assume the search space $\X$ to be finite and that,
w.l.o.g., all objective shall be maximized. We are therefore
interested in solving
\begin{align}
\label{MOP}
\text{max} \left\{ f(x) = (f_1(x), \ldots, f_m(x))^\T \mid x \in \X
\right\} 
\end{align}
Here, maximization is understood with respect to the product order
$\ge$ on $\Real^m$, i.e., for any pair $(y,y') \in \Real^m \times
\Real^m$, $y \ge y'$ if and only if $y_i \ge y'_i$ for all $i \in \{1,
\ldots, m\}$. Hence $(\Y,\ge)$, where $\Y = f(\X)$ is called the
objective space, is a partially ordered set. This gives rise to the
following order relation on the search space.  

\begin{definition}[Pareto dominance]
For any pair $(x,x') \in \X \times \X$, $x$ is said to weakly dominate
$x'$, denoted as $x \succeq x'$, if and only if $f(x) \ge f(x')$. $x$
said to dominate $x'$, denoted as $x \succ x'$, if and only if $x
\succeq x'$ and $x' \not\succeq  x$.
\end{definition}

Note that $(\X,\succeq)$ is a preordered set, while $(\X,\succ)$
is a strictly partially ordered set. 
A subset $\X' \subseteq \X$ is called \emph{independent} with respect
to $\succeq$ if for all $(x,x') \in \X' \times \X'$ where $x \not= x'$
it holds that $x \not\succeq x'$. Let the set of all such independent
subsets of $\X$ be denoted by $\hat{\X}$, i.e.,
\[
\hat{\X} := \{ \X' \subseteq \X \mid \X' \text{ is independent with
  respect to } \succeq \}.
\]
The set of minimal elements of the objective space,  
\[
\Y^* := \max(\Y,\ge) = \{ y \in \Y \mid \not\exists y' \not= y 
\text{ with } y' \ge y\}, 
\] 
is called the \emph{Pareto front}, and its preimage $\X^* =
f^{-1}(\Y^*)$ is called the \emph{Pareto set}.

Ideally, when solving (\ref{MOP}), one is interested in determining
the Pareto front $\Y^*$, together with an independent set $\X'$ that
should \emph{cover} $\Y^*$, i.e, $f(\X') = \Y^*$. This means we are
usually not interested in obtaining more than one preimage for each
element of the Pareto front. 

In many instances, the size of Pareto front might be immense, so we
are interested in approximations. Here, our goal is for find some
subset of $\Xh$ of a given maximum size $k$ that approximates $Y^*$
well in the following sense.

\begin{definition}[$\e$-dominance]
\label{def:e}
Let $\e \in \Real^m, \e \geq 0$. For any pair $(x,x') \in \X \times
\X$, $x$ is said to $\e$-dominate $x'$, denoted as $x \se x'$,
if and only if $f(x) + \e \ge f(x')$.
\end{definition}

A set $\A \in \hat{\X}$ is called an $\e$-approximate Pareto set if
for all $x \in \X$ there is an $x' \in \A$ with $x' \se x$. An
$\e$-approximate Pareto set that is a subset of $X^*$ is called an
$\e$-Pareto set. Thus, a reasonable task is to find an $\e$-Pareto set
of some given maximum cardinality $k$. Note that for the special case
of $\e=0$, the notions of $\e$-approximate Pareto set, $\e$-Pareto
set, and covering independent set are equivalent.

We investigate in which sense simple random search is able to find
an $\e$-Pareto set of cardinality $k$ of the multiobjective
optimization problem (\ref{MOP}), i.e., whether its solution set
stochastically converges to such a set in the limit. 
For this, we consider Algorithm~\ref{RS}, which is pure random search
where the set $\A_t$ represents its \emph{archive} of solutions,
stored in a memory (array) of size at most $k$.

\begin{algorithm}[t]
  \caption{Simple random search algorithm}
  \label{RS}
  \begin{algorithmic}[1]
    \STATE $t \leftarrow 0$
    \STATE $\A_t \leftarrow \emptyset$
    \LOOP
    \STATE $t \leftarrow t + 1$
    \STATE Draw $x_t$ from $\X$ uniformly at random
    \STATE $\A_t \leftarrow \up(\A_{t-1},x_t)$
    \ENDLOOP
  \end{algorithmic}
\end{algorithm}

\section{Related work and previous results}

Many different notions of Pareto set approximations have been proposed
in the literature, see for instance the survey on concepts of
$\epsilon$-efficiency in \cite{HP1994a}. However, many of them deal with
infinite sets and are therefore not practical as a solution concept for our
purpose of producing and maintaining a representative subset of a
fixed given size. The use of discrete $\epsilon$-approximations of the
Pareto set was proposed almost simultaneously by various authors
\cite{hansen_1980_bicriterion,EP1987a,Reuter1990a,RF1990a}. The
general idea is that each Pareto-optimal point is approximately
dominated by some point of the approximation set. The $\e$-dominance
given in Definition~\ref{def:e} is a typical instance of this
approach, while it is also common to use a relative deviation instead
of the absolute deviation employed here.  

As relative deviation is essentially equivalent to absolute deviation
on a logarithmically scaled objective space, this choice should not
affect the convergence results obtained but rather depend on
the actual application problem at hand. The nice property of relative
deviation is that it allows to prove that, under very mild
assumptions, there is always an $\e$-Pareto set whose size is
polynomial in the input length
\cite{PY2000,erlebach_2002_approximating}. Further approximation 
results for particular combinatorial multiobjective optimization
problems are given in \cite{Ehrg2000a}, where the question was how well
a single solution can approximate the whole Pareto set, which is a
special case of our question restricted to $k=1$ and with focus on 
deterministic algorithms.

Despite the existence of suitable approximation concepts,
investigations on the \emph{convergence} of particular algorithms
towards such approximation sets, that is, their ability to obtain a
suitable Pareto set approximation in the limit, have remained
rare. In \cite{Rudolph1998c,RA2000a} the stochastic search procedure
proposed by earlier by \cite{peschel_1977_use} was analyzed and proved
to converge to an $\e$-Pareto set with $\e=0$ in case of a finite
search space. Obviously, the solution set maintained by this algorithm
might in the worst case grow as large as the Pareto set $\X^*$
itself. Thus, a different version with bounded memory of at most $k$
elements was proposed and shown to converge to some subset of $\X^*$ of
size at most $k$, but no guarantee about the approximation quality
could be given. Similar results were obtained by \cite{Hanne1999} for
continuous search spaces. 

One option to control the approximation quality under size
restrictions is to define a \emph{quality indicator} which maps each
possible solution set to a real value that can then be used to decide
on the inclusion of a new search point. Several algorithms have been
proposed that implement this concept \cite{zk2004a,beume_2007_sms}.
In case that such a quality indicator  
fulfils certain monotonicity conditions, it can be used as a potential
function in the convergence analysis. As shown in
\cite{knowles_2002_local,knowles_2003_properties}, this entails 
convergence to a subset of the Pareto set as a local optimum
of the quality indicator, but it remained open how such a local
optimum relates to a guarantee on the approximation quality
$\e$. \cite{knowles_2003_properties} also analyzed an adaptive grid 
archiving method proposed in \cite{knowles_2000_approximating} and
proved that after finite time, even though the solution set itself
might permanently oscillate, it will always represent an
$\epsilon$-approximation whose approximation quality depends on the
granularity of the adaptive grid and on the number of allowed
solutions. The results depend on the additional assumption that the
grid boundaries converge after finite time, which is fulfilled in
certain special cases. 

In \cite{LTDZ2002b}, two archiving algorithms were proposed
that provably maintain a finite-size $\e$-approximation of all points
ever generated during the search process. As an immediate corollary,
these archiving strategies ensure convergence to a Pareto set
approximation of given quality $\e$. While the desired $\e$ is an
input to the algorithm that can be specified beforehand, the resulting
size of the approximation can be bounded as a function of $\e$ and the
ranges of the objective space $\Y$. In \cite{knowles_2004_bounded} it
was shown that in practise these bounds are sometimes not tight
enough, which is a 
particular disadvantage when applied in a scenario where the maximum
archive size $k$ has to be specified beforehand. 
If information of the objective space ranges is available, the bounds
can be used to find a valid value for $\e$, but this choice
often turns out overly conservative so that far less solutions are
attained than would be possible with a memory of size $k$. In case
where the objective ranges are unknown, a mechanism proposed in
\cite{LTDZ2002b} can be used to systematically increase the
value of $\e$ without losing the convergence properties, but this
suffers from the same drawback of being overly conservative. Thus, it
has remained open until now whether working with fixed size Pareto
set approximations can guarantee convergence in the limit for
arbitrary multiobjective optimization problems on finite search spaces,
and at the same time guarantee a certain approximation quality.

In this paper we settle this question positively by presenting an 
archiving scheme that enables Algorithm~\ref{RS} to produce a sequence
of solution sets that converges with probability 
one to an $\e$-Pareto set of a certain quality. The algorithm
represents a combination of the two complementary algorithms discussed
above \cite{knowles_2003_properties,LTDZ2002b}, thus combining the
advantages of 
both: it can be seen as a variant of the adaptive grid archiving
method where a multilevel, fixed grid is used instead of an adaptive
grid with moving boundaries, which is crucial to obtain
convergence. It can also be seen as a proper implementation of the
adaptation mechanism for the $\e$ values proposed in \cite{LTDZ2002b},
which is crucial to limit the size of the solution set to at most $k$,
but which is also able to reduce the value of $\e$ whenever
possible. 
Finally, the algorithm can be seen as selection using a particular
quality indicator \cite{zk2004a}, a notion that will be defined more
precisely later on. However, instead of 
having to compute the actual indicator values, which might be
computationally cumbersome, this indicator will only be
used as a potential function in the analysis of the algorithm and
never has to be computed. The actual comparison will be defined using
very simple local rules that -- and this is crucial -- are in
accordance with the quality indicator, which will be established via
order homomorphisms. 

\section{Stochastic convergence analysis}

The sequences of archives $\{\A_t, t \in \Nat_0\}$ generated by
Algorithm~\ref{RS} are realizations of a discrete-time stochastic
process defined on a probability space $(\Omega,\F,\Prob)$. The sample
space $\Omega$ can be defined as the infinite product
$\hat{\X}^\infty$ with the $\sigma$-algebra $\F = 2^\Omega$ being its
power set. The probability measure $\Prob$ is defined implicitly by
Algorithm~\ref{RS}. The stochastic process is then the sequence of
random variables $\{A_t, t \in \Nat_0\}$, where 
\[
A_t : \Omega \lra \hat{\X}, \quad A_t(\omega) = \omega_t
\]
and $\A_t = A_t(\omega)$. 

From Algorithm~\ref{RS} is it clear that $A_{t+1}$ only depends on
$A_t$, so the process is a homogenous finite Markov chain with state space
$\hat{\X}$. Due to line~5 of the algorithm, all transition
probabilities that change the state are equal to $p = 1/|\X|$. 

Let the transition graph of the Markov chain be denoted by
$G=(\hat{\X},E)$. Clearly, its arcs $E \subseteq \hat{\X} \times
\hat{\X}$ are determined by the $\up$ function given
in Algorithm~\ref{up} as
\[
E = \{ (\A,\A') \in \hat{\X} \times \hat{\X} \mid \exists x \in \X :
\A' = \up(\A,x) \}.
\]

Our goal is to show that this transition graph, when ignoring loops,
forms a directed acyclic graph, which immediately implies that
absorption will take place with probability one in finite time. 
We then proceed to show that in all absorbing states the archives are
$\e$-Pareto sets, and finally give some guarantee of the approximation
quality $\e$ obtained. 

Instead of working on $E$ directly, however, we define a potential function
$I$, according to which the set of possible archives $\hat{\X}$ can be
linearly ordered. For this, some auxiliary notation is needed.

\begin{definition}[box index vector]
The box index vector of a vector $y \in \Real^m$ at level $b \in \Zed$ is
given by the value of the function
\[
\beta^{(b)}: \Real^m \lra \Zed^m, \quad
\beta^{(b)}(y) = \left(r^{(b)}(y_1),\ldots,r^{(b)}(y_m) \right)^\T
\]
where
\[
r^{(b)}: \Real \lra \Zed, \quad
r^{(b)}(z) = \left\lfloor z \cdot 2^{-b} +  \frac{1}{2} \right\rfloor
\]
\end{definition}

\begin{definition}[box-dominance]
Let $b \in \Zed$. For any pair $(x,x') \in \X \times \X$, $x$ is said
to weakly box-dominate $x'$ at level $b$, denoted as $x \be x'$, if
and only if $\beta^{(b)}(f(x)) \ge \beta^{(b)}(f(x'))$. 
$x$ is said to be box-equal to $x'$, denoted as $x \sim_b x'$, if $x
\be x'$ and $x' \be x$. If $x \be x'$ and $x \not\sim_b x'$ then $x$
is said to box-dominate $x'$ at level $b$, denoted as $x \rhd_b x'$. 
\end{definition}

Note that the relations $\be$ form a family of order extensions of
the dominance relation $\succeq$. The accompanying equivalence relations
$\sim_b$ can be seen as a successive coarse-graining of approximate
indifference between solutions.

\begin{lemma}
\label{lem:order-extension}
If $x \be x'$ then $ x \unrhd_c x'$ for all $c \ge b$.
\end{lemma}

\begin{proof}
We show for all components $i \in \{1,\ldots,m\}$ that if
$r^{(b)}(f_i(x)) \ge r^{(b)}(f_i(x))$ then $r^{(b+1)}(f_i(x)) \ge
r^{(b+1)}(f_i(x))$. The lemma follows then by induction. Let $d :=
f_i(x) \cdot 2^{-b} + 1/2$ and $d' := f_i(x') \cdot 2^{-b} + 1/2$. From
the premise we can express $d = p + 1/2 + g$ and $d' = p + 1/2 + 1 -
h$ for some $p \in \Zed$, $g \ge 0$ and $h > 0$. If $p = 2k$ for some
$k \in \Zed$ then $r^{(b+1)}(d) \ge \lfloor 2k/2 + 1/2  \rfloor = k \ge
\lfloor 2k/2 + 1 - d \rfloor = r^{(b+1)}(d')$. If $p=2k+1$ then
$r^{(b+1)}(d) \ge \lfloor (2k+1)/2 + 1/2  \rfloor = k+1 \ge \lfloor
(2k+1)/2 + 1 - d \rfloor = r^{(b+1)}(d')$. 
\end{proof} 

\begin{definition}[potential function]
Let 
\[
I: \hat{\X} \lra \Real, \quad
I(\A) := \sum_{i = -\bar{b}}^{\infty} |\B_{-i}(\A)| \cdot (|\X|+1)^{-i}
\]
where
\[
\B_b(\A) = \{ x \in \X \mid \exists a \in \A : a \be x \}
\]
and
\[
\bar{b} = \min \{ b \in \Zed \mid \forall (x,x') \in \X \times \X : 
x \sim_b x'\}
\]
\end{definition}
The power series defining $I$ converges as the $\B_b$ are subsets of
$\X$, which is finite. Moreover, $\bar{b}$ exists since it is possible
to enclose the whole objective space $f(\X)$ by one box by choosing
$b$ large enough. 

\begin{algorithm}[t]
  \caption{$\up(\A,x)$}
  \label{up}
  \begin{algorithmic}[1]
    \STATE $\A' \leftarrow \A \cup \{x\}$
    \IF {$\exists a \in \A : a \succeq x$}
    \STATE \textbf{return} $\A$
    \ELSIF {$\max(\A', \succ) \leq k$}
    \STATE \textbf{return} $\max(\A', \succ)$
    \ELSE
    \STATE $\beta \leftarrow \min\{b \in \Zed \mid \exists (a,a') \in \A'
    \times \A' : a \be a' \text{ and } a \not= a'\}$
    \STATE $\A'' \leftarrow \{ a \in \A' \mid \exists a' \in \A' : 
    a' \be a \text{ and } a' \not= a\}$ 
    \IF {$x \in \A''$}
    \STATE \textbf{return} $\A$
    \ELSE
    \STATE Draw $a$ uniformly at random from $\A''$
    \STATE \textbf{return} $\A \setminus \{a\} \cup \{x\}$
    \ENDIF
    \ENDIF
  \end{algorithmic}
\end{algorithm}

The dominance relation on solutions can be used to define a natural
preference relation on the set of independent sets $\hat{\X}$.

\begin{definition}[dominance of independent sets]
Let $(\A,\A') \in \hat{\A} \times \hat{\A}$. The set $\A$ is said to
weakly dominate $\A'$, denoted as $\A \sqsupseteq \A'$, if $\max(f(\A
\cup \A'),\ge) = \max(f(\A),\ge)$. 
\end{definition}

\begin{lemma}
\label{lem:qi}
$I$ is an order homomorphism of $(\hat{\A},\sqsupseteq)$ into
$(\Real,\ge)$, i.e., if $\A \sqsupseteq \A'$ then $I(\A) \ge
I(\A')$. $I$ is also an order homomorphism of $(\hat{\A},\sqsupset)$
into  $(\Real,>)$, i.e., if $\A \sqsupset \A'$ then $I(\A) > I(\A')$. 
\end{lemma}

\begin{proof}
If $\A \sqsupseteq \A'$ then the coefficients $|\B_i(\A)|$ in the power
series of $I(\A)$ are uniformly not less than $|\B_i(\A')|$ because
$\B_i(\A') \subseteq \B_i(\A)$ for all $i$. If additionally $\A'
\not\sqsupseteq \A$ then there exists an $a \in \A$ and a $b \in \Zed$
such that there is no $a' \in \A'$ with $a' \be a$. Hence, $\B_b(\A')
\subset \B_b(\A)$, which implies that $I(\A) > I(\A')$.
\end{proof}

The potential function $I$ can bee seen as a quality indicator for
independent sets. As an immediate corollary of Lemma~\ref{lem:qi}, the
$>$-relation on $I(\hat{\A})$ represents a comparison method that is 
$\not\sqsupseteq$-compatible and $\sqsupset$-complete in the
terminology of \cite{ZTLFF2003a}. 

\begin{lemma}
\label{lem:dag}
If $(\A_t,\A_{t+1}) \in E$ and $\A_t \not= \A_{t+1}$ then $I(\A_{t+1})
> I(\A_t)$.
\end{lemma}
 
\begin{proof}
If $(\A_t,\A_{t+1}) \in E$ and $\A_t \not= \A_{t+1}$ then $x \in
\A_{t+1} = \up(\A_t,x)$. There are two cases that can cause
$x$ to be included into the new archive (termination in line~5 or
line~13). For termination to occur in line~5, $\max(\A', \succ) \leq 
k$ has to hold (line 4), and furthermore $x \in \max(\A_{t+1},\succ)$
by contradiction to the condition in line~2.  
Since $\A_t = \max(\A, \succ) \not\ni x$ it follows that $\A_{t+1} =
\max(\A', \succ) \sqsupset \A_t$ and thus, by Lemma~\ref{lem:qi},
$I(\A_{t+1}) > I(\A_t)$. When termination occurs in line~13,
$\max(\A',\unrhd_\beta)$ contains $x$ but not $a$. Hence
$\B_\beta(\A_{t+1}) \supset \B_\beta(\A_t)$, which implies, by
Lemma~\ref{lem:order-extension}, that $\B_b(\A_{t+1}) \supseteq
\B_b(\A_t)$ for all $b \ge \beta$. This implies $I(\A_{t+1}) >
I(\A_t)$, which completes the proof.
\end{proof}

\begin{theorem}
\label{the:absorbing}
The Markov chain $\{A_t, t \in \Nat_0\}$ is absorbing.
\end{theorem}

\begin{proof}
Due to Lemma~\ref{lem:dag}, the mapping $I$ is an order homomorphism
of $G$ into a strict linear order. This implies that the transitive
closure of $G$ is a partially ordered set, and hence $G$ is a directed
acyclic graph. As any Markov chain on an directed acyclic graph is
absorbing, the claim follows.
\end{proof}

\begin{theorem}
\label{the:pareto}
For any absorbing state $\A \in \max(G)$ it holds that $\A
\subseteq \X^*$.
\end{theorem}

\begin{proof}
Assume that $\A \not\subseteq \X^*$, so there is some $a \in \A$ that
is not in $\X^*$. Then there exists some $x \in \X$ with $x \succ
a$. Let $\tilde{\A} = \up(\A,x)$. Since $x \succ A$, the condition
$\max(\A',\succ) \le k$ in line~4 holds because $a \not\in
\max(\A',\succ)$. Thus, $\tilde{\A} = \max(\A',\succ) \not= \A$, which
is a contradiction to the assumption.
\end{proof}

The next theorem shows that any absorbing state box-dominates the
Pareto set at the lowest level possible with the least number of
solutions necessary, while distributing the remaining solutions
with maximum entropy over the nondominated boxes at the next lower
level. 
\begin{theorem}
\label{the:b}
Let
\begin{align}
\label{eq:delta}
\delta & = \min \{ b \in \Zed \mid \exists \A \subseteq \X^* \text{ with
  } |\A| \le k \text{ and } \B_b(\A) = \X \},
\end{align}
denote the smallest level $b$ on which it is possible to box-dominate
the Pareto set with at most $k$ solutions,
\begin{align}
s & = \min \{ |\A| \, : \, \B_\delta(\A) = \X \},
\end{align}
denote the minimum size of such a set, and let
\begin{align}
\E & = \{ \X' \subseteq \X^* \mid \forall (x,x') \in \X' \times \X' : x
  \sim_{\delta-1} x' \}
\end{align}
denote a partitioning of the Pareto set into the boxes of the next
smaller level.
Then for any absorbing state $\A \in \max(G)$ it holds that
$\B_\delta(\A) = \X$ and that
\begin{align}
\label{ks}
|\{ \X' \in \E : \X' \cap \A \not= \emptyset \}| = \min\{k,|\E|\}.
\end{align}
\end{theorem}

\begin{proof}
Assume that there is some $x \in \X^*$ with $x \not\in
\B_\delta(\A)$. If $|\A| < k$ than $\max(\A', \succ) \le k$ in line~4,
in which case $\A$ cannot be absorbing. 
If $|\A| = k$ then $\A''$ cannot be empty, as this would contradict
the definition of $\beta$ in line~7. If $x \not\in \A''$ then $x$ will
enter $\A$, so $\A$ cannot be absorbing. 
If $x \in \A''$ then $\beta > \delta$ due to
Lemma~\ref{lem:order-extension}, since $x \not\in \B_\delta(\A)$ but
$x \in \B_\beta(\A)$ by assumption.
The definition of $\beta$ in line 7 now implies that $\A'$ is an
independent set of cardinality $k+1$ with respect to $\unrhd_\delta$,
hence $\A'$ serves as a witness for the fact that there is no
$\A$ with $|\A| \le k$ and $\B_b(\A) = \X$, which is a contradiction to
(\ref{eq:delta}) and completes the proof that $\B_\delta(\A) = \X$. 
To prove the second part of the proposition, assume that
$\B_\delta(\A) = \X$ and $d := |\{ \X' \in \E : \X' \cap \A \not=
\emptyset \}| < \min\{k,|\E|\}$. Hence there is some $x \in X' \in
\E$ with $\X' \cap \A = \emptyset$. If $|\A| < k$ than $x$ will be
accepted as there is no $a \in \A$ with $a \succeq x$, so $\A$ cannot
be absorbing. If $|\A| = k$ then $\beta < \delta$ as otherwise $d =
k$. Now, as $\not\exists a \in \A$ with $a \unrhd_{\delta-1} x$ and $\beta
\le \delta-1$, it follows with Lemma~\ref{lem:order-extension} that $x
\not\in \A''$, thus $x$ will be accepted, so $\A$ cannot be
absorbing.
\end{proof}

We have collected now all necessary ingredients to show the main
result. 

\begin{theorem}
The sequence $\{A_t, t \in \Nat_0\}$ converges with probability one to
some $\e$-Pareto set with $\e = 2^\delta$, with $\delta$ defined as in
Theorem~\ref{the:b}.  
\end{theorem}

\begin{proof}
As a corollary to Theorems~\ref{the:absorbing} and \ref{the:pareto},
$A_t$ converges with probability one to some subset of the Pareto set.
As a corollary of Theorem~\ref{the:b}, for any absorbing state $\A$
there exists for all $x \in X^*$ some $a \in \A$ such that $a
\unrhd_\delta x$ and, hence, $a \succ_\e x$ for all $\e \ge 2^\delta$.
\end{proof}

\section{Conclusions}

In this paper, the first convergence result for random search
algorithms to a subset of the Pareto set of given maximum size $k$ and
bounds on the approximation quality $\e$ was given. The convergence
was enabled by a new selection scheme, given als Algorithm~\ref{up},
that compares the new candidate solution to the current archive using
a multi-level grid.

In many parts, the
assumption of a finite search space was used. Even though this is a
reasonable assumption for any implementation in computer arithmetic
with finite precision, an extension to the continuous case would be
desirable. Even though it might be a justified assuption that the
results can be extended, recent experience
\cite{schuetze_2007_convergence} has shown that this might involve
considerable effort. 

\bibliographystyle{plain}
\bibliography{pareto,knowles,schuetze,erlebach}

\end{document}